\documentclass{amsart}

\usepackage{fancyhdr}
\usepackage{lastpage}
\usepackage{stmaryrd}

\newcommand{\bdis}{\begin{displaymath}}
\newcommand{\edis}{\end{displaymath}}
\newcommand{\be}{\begin{equation}}
\newcommand{\ee}{\end{equation}}

\newcommand{\mcal}{\mathcal}

\newtheorem{theorem}{Theorem}

\theoremstyle{definition}

\newtheorem{cor}[]{Corollary}

\theoremstyle{remark}
\newtheorem{remark}[]{Remark}

\newtheorem*{mydef3}{{\bf Problem}}

\numberwithin{equation}{section}

\begin{document}

\title{Jacob's ladders and the quantization of the Hardy-Littlewood integral}

\author{Jan Moser}

\address{Department of Mathematical Analysis and Numerical Mathematics, Comenius University, Mlynska Dolina M105, 842 48 Bratislava, SLOVAKIA}

\email{jan.mozer@fmph.uniba.sk}

\keywords{Riemann zeta-function}

\begin{abstract}
We use Jacob's ladders to solve the fine problem how to divide of the Hardy-Littlewood integral to equal parts,
for example  of magnitude $h=6.6\times 10^{-27}$ (the numerical value of elementary Planck quantum). The result of the paper cannot be obtained in known
theories of Balasubramanian, Heath-Brown and Ivic.
\end{abstract}

\maketitle

\section{The problem of dividing on equal parts}

\subsection{}

Let us remind the following facts. Titchmarsh-Kober-Atkinson (TKA) formula
\be \label{form1.1}
\int_0^\infty Z^2(t)e^{-2\delta t}{\rm d}t=\frac{c-\ln (4\pi\delta)}{2\sin \delta}+\sum_{n=1}^N c_n\delta^n+\mcal{O}(\delta^{N+1})
\ee
(see \cite{10}, p. 141) remained as an isolated result for a period of 56 years until we have discovered the nonlinear integral equation
\be \label{form1.2}
\int_0^{\mu[x(T)]}Z^2(t)e^{-\frac{2}{x(T)}t}{\rm d}t=\int_0^TZ^2(t){\rm d}t
\ee
(see \cite{6}) in which the essence of the TKA formula is encoded. Namely, we have shown in \cite{6} that the following almost exact formula for the
Hardy-Littlewood integral takes place
\be \label{form1.3}
\int_0^TZ^2(t){\rm d}t=\frac{\varphi(T)}{2}\ln\frac{\varphi(T)}{2}+(c-2\pi)\frac{\varphi(T)}{2}+c_0+\mcal{O}\left(\frac{\ln T}{T}\right),
\ee
where $\varphi(T)$ is the Jacob's ladder (a solution of the nonlinear integral equation (\ref{form1.2})).

\begin{remark}
Our formula (\ref{form1.3}) for the Hardy-Littlewood integral has been obtained after the time period of 90 years since this integral appeared in 1918
(see \cite{3}, pp. 122, 151-156).

\end{remark}

\begin{remark}
Let us remind that
\begin{itemize}
\item[(A)] The Good's $\Omega$-theorem (see \cite{2}) implies for the Balasubramanian formula
\be \label{form1.4}
\int_0^TZ^2(t){\rm d}t=T\ln T+(2c-1-\ln 2\pi)T+R(T),\ R(T)=\mcal{O}(T^{1/3+\epsilon})
\ee
(see \cite{1}) that
\bdis
\limsup_{T\to\infty}|R(T)|=+\infty .
\edis
\item[(B)] The error term in (\ref{form1.3}) tends to zero as $T$ goes to infinity, namely
\bdis
\lim_{T\to\infty}r(T)=0,\ r(T)=\mcal{O}\left(\frac{\ln T}{T}\right),
\edis
i.e. our formula is almost exact (see \cite{6}).
\end{itemize}
\end{remark}

\subsection{}

In this paper I consider the problem concerning the solid of revolution corresponding to the graph of the function $Z(t),\ t\in [T_0,+\infty)$, where
$0<T_0$ is a sufficiently big number.

\begin{mydef3}

To divide this solid of revolution on parts of equal volumes.

\end{mydef3}

We obtain, for example, from our formula (\ref{form1.3}) that there exists a sequence $\{ \hat{T}_\nu\}_{\nu=\nu_0}^\infty$, for which
\be\label{form1.5}
\pi\int_{\hat{T}_\nu}^{\hat{T}_{\nu+1}}Z^2(t){\rm d}t=6.6\times 10^{-27} ,
\ee
\bdis
\hat{T}_{\nu+1}-\hat{T}_{\nu}\sim \frac{6.6\times 10^{-27}}{\pi\ln \hat{T}_\nu\tan[\alpha(\hat{T}_\nu,\hat{T}_{\nu+1})]},\ \nu\to\infty ,
\edis
where $h=6.6\times 10^{-27}\text{erg}\cdot \text{sec}$ is elementary Planck quantum.

\begin{remark}

It us quite evident that the quantization rule (\ref{form1.5}) cannot be obtained by methods of Balasubramanian, Heath-Brown and Ivic (see, for example,
\cite{4}).

\end{remark}

This paper is a continuation of the series of papers \cite{6}-\cite{8}.

\section{Main result}

The following theorem holds true

\begin{theorem}
Let $0<\delta<\Delta$ where $\delta$ is an arbitrarily small and $\Delta$ is an arbitrarily big number. Then for every $\omega\in [\delta,\Delta]$
and every Jacob's ladder $\varphi(T)$ there is the sequence
\bdis
\{T_\nu(\omega,\varphi)\}_{\nu=\nu_0}^\infty,\quad T_\nu(\omega,\varphi)=T_\nu(\omega)
\edis
for which
\be\label{form2.1}
\int_{T_\nu(\omega)}^{T_{\nu+1}(\omega)}Z^2(t){\rm d}t=\omega,
\ee
\be\label{form2.2}
T_{\nu+1}(\omega)-T_\nu(\omega)=\frac{\omega+\mcal{O}\left(\frac{\ln T_\nu}{T_\nu}\right)}
{\left( \ln\frac{\varphi(T_\nu)}{2}-a\right)\tan[\alpha(T_\nu,T_{\nu+1})]} ,
\ee
where $0<\nu_0(\omega,\varphi)$ is a sufficiently big number, $a=\ln 2\pi-1-c$ and $\alpha=\alpha(T_\nu,T_{\nu+1})$ is he angle of the chord binding
the points
\bdis
\left[ T_\nu,\frac{1}{2}\varphi[T_\nu(\omega)]\right],\quad
\left[ T_{\nu+1},\frac{1}{2}\varphi[T_{\nu+1}(\omega)]\right]
\edis
of the curve $y=\frac{1}{2}\varphi(T)$.
\end{theorem}

We obtain from our Theorem 1

\begin{cor}
\begin{eqnarray} \label{form2.3}
& & \int_{T_\nu(\omega)}^{T_{\nu+1}(\omega)}Z^2(t){\rm d}t= \\
& &
[T_{\nu+1}(\omega)-T_\nu(\omega)]\ln\left( e^{-a}\frac{\varphi(T_\nu)}{2}\right)\tan[\alpha(T_\nu,T_{\nu+1})]+
\mcal{O}\left(\frac{\ln T_\nu}{T_\nu}\right). \nonumber
\end{eqnarray}
\end{cor}

Let us remind the formula (see \cite{7}, (2.1))
\be\label{form2.4}
\int_T^{T+U}Z^2(t){\rm d}t=U\ln\left( e^{-a}\frac{\varphi(T)}{2}\right)\tan[\alpha(T,U)]+\mcal{O}\left(\frac{1}{T^{1/3-4\epsilon}}\right).
\ee

\begin{remark}

In the case of collection of sequences $\{ T_\nu(\omega)\}$ we obtain the essential improvement
\bdis
\mcal{O}\left(\frac{1}{T^{1/3-4\epsilon}}\right)\rightarrow \mcal{O}\left(\frac{\ln T_\nu}{T_\nu}\right)
\edis
(see (\ref{form2.3}) of the remainder term in the formula (\ref{form2.4}).

\end{remark}

Since
\bdis
\sum_{k=1}^N \int_{T_{\nu+k-1}(\omega)}^{T_{\nu+k}(\omega)}Z^2(t){\rm d}t=N\omega
\edis
then, in the case
\bdis
T_{\nu+N_0}(\omega)-T_{\nu}(\omega)\sim U_0=T^{1/3+2\epsilon} ,
\edis
we have (see (\ref{form1.4}))
\bdis
N_0\sim \frac{1}{\omega}U_0\ln T_\nu(\omega),\ \frac{U_0}{N_0}\sim \frac{\omega}{\ln T_\nu(\omega)} .
\edis

Then we obtain by our Theorem 1

\begin{cor}

The following asymptotic formula takes place
\be\label{form2.5}
\frac{1}{N_0}\sum_{k=1}^{N_0}\left\{ T_{\nu+k}(\omega)-T_{\nu+k-1}(\omega)\right\}\sim \frac{\omega}{\ln T_\nu(\omega)}
\ee
for arithmetic mean values of $T_{\nu+k}(\omega)-T_{\nu+k-1}(\omega)$, $k=1,2,\dots ,N_0$.

\end{cor}

\section{On transformation of the sequence $\{ T_\nu(\omega)\}$ preserving the quantization of the Hardy-Littlewood integral}

\subsection{}

If there is a sufficiently big natural number $\bar{\nu}$ for which $\omega\bar{\nu}=T_0$ is fulfilled then from our Theorem 1 the resolution of our
Problem follows. The complete resolution follows from the next theorem.

\begin{theorem}
For every $\omega\in [\delta,\Delta],\ \tau\in [0,\omega)$ and every Jacob's ladder $\varphi(T)$ there is the collection of sequences
\bdis
\{ T_\nu(\omega,\tau;\varphi)\}_{\nu=\nu_0}^\infty, \quad
T_\nu(\omega,\tau;\varphi)=T_\nu(\omega,\tau),\ T_\nu(\omega,0)=T_\nu(\omega)
\edis
for which
\be \label{form3.1}
\int_{T_\nu(\omega,\tau)}^{T_{\nu+1}(\omega,\tau)}Z^2(t){\rm d}t=\omega ,
\ee
\be \label{form3.2}
T_{\nu+1}(\omega,\tau)-T_{\nu}(\omega,\tau)=\frac{\omega+\mcal{O}\left(\frac{\ln T_\nu}{T_\nu}\right)}
{\left(\ln\frac{\varphi(T_\nu)}{2}-a\right)\tan[\alpha(T_\nu,T_{\nu+1})]} ,
\ee
where $0<\nu_0(\omega,\tau;\varphi)$ is a sufficiently big number and $\alpha=\alpha(T_\nu,T_{\nu+1})$ is the angle of the chord binding the points
\bdis
\left[ T_\nu(\omega,\tau),\frac{1}{2}\varphi(T_\nu(\omega,\tau))\right],\quad
\left[ T_{\nu+1}(\omega,\tau),\frac{1}{2}\varphi(T_{\nu+1}(\omega,\tau))\right]
\edis
of the curve $y=\frac{1}{2}\varphi(T)$.
\end{theorem}

\begin{remark}

From (\ref{form3.1}), (\ref{form3.2}) we get full analogies of Corollaries  1., 2. and Remark 4.

\end{remark}

\begin{remark}

The quantization (\ref{form1.5}) follows from the choice $\{ T_\nu(\omega,\tau)\}$ with $\omega=\frac{h}{\pi},\ \bar{\nu}\omega+\tau=T_0$.

\end{remark}

\subsection{}

Let us remind that for Gram's sequence $\{ t_\nu\}$ we have (see \cite{9})
\be \label{form3.3}
t_{\nu+1}-t_\nu \sim \frac{2\pi}{\ln t_\nu} ,
\ee
and for the collection of sequences $\{ \bar{t}_\nu(\bar{\tau})\},\ \bar{\tau}\in [-\pi,\pi]$ defined in our paper \cite{5} we have the analogue of
(\ref{form3.3})
\bdis
\bar{t}_{\nu+1}(\bar{\tau})-\bar{t}_{\nu}(\bar{\tau})\sim \frac{2\pi}{\ln t_\nu(\bar{\tau})} .
\edis

\begin{remark}

Under the transformations
\bdis
t_\nu\rightarrow T_\nu(\omega);\ T_\nu(\omega,\tau)
\edis
the individual property (3.3) is transformed to an analogous property of arithmetic means (see (\ref{form2.5}) and Remark 5).

\end{remark}

\section{Proof of Theorems 1., 2.}

\subsection{}

By (\ref{form1.3}) we have
\be \label{form4.1}
\int_0^T Z^2(t){\rm d}t=F[\varphi(T)]+\mcal{O}\left(\frac{\ln T}{T}\right),\ T\geq T^{(1)}[\varphi] ,
\ee
\be \label{form4.2}
F(y)=\frac{y}{2}\ln\frac{y}{2}+(c-\ln 2\pi)\frac{y}{2}+c_0,\
F'(y)=\frac{1}{2}\ln\frac{y}{2}-\frac{a}{2},\ F''(y)=\frac{1}{2y} .
\ee
Since the continuous function
\bdis
F[\varphi(T)]+\mcal{O}\left(\frac{\ln T}{T}\right),\ T\geq T^{(1)}
\edis
is increasing, there is a root $T_\nu(\omega,\varphi)$ of the equation
\bdis
F[\varphi(T)]+\mcal{O}\left(\frac{\ln T}{T}\right)=\omega \nu,\quad \nu\geq \nu_0,
\edis
where $\nu_0=\nu_0(\omega,\varphi)$ is a sufficiently big number. Thus, the sequence $\{ T_\nu(\omega;\varphi)\}_{\nu=\nu_0}^\infty$ is
constructed by equation
\be \label{form4.3}
F[\varphi(T_\nu)]+\mcal{O}\left(\frac{\ln T_\nu}{T_\nu}\right)=\omega \nu .
\ee
From (\ref{form4.1}) by (\ref{form4.3}) we obtain
\bdis
\int_0^{T_\nu(\omega;\varphi)}Z^2(t){\rm d}t=\omega\nu \ \Rightarrow\ \int_{T_\nu(\omega;\varphi)}^{T_{\nu+1}(\omega;\varphi)}Z^2(t){\rm d}t=\omega,
\edis
i.e. (\ref{form2.1}).

\subsection{}

We have by (\ref{form4.2}), (\ref{form4.3})
\begin{eqnarray*}
& &
\omega+\mcal{O}\left(\frac{\ln T_\nu}{T_\nu}\right)=F[\varphi(T_{\nu+1})]-F[\varphi(T_\nu)]= \\
& &
\left( \frac{1}{2}\ln\frac{d}{2}-\frac{a}{2}\right)\left[ \varphi(T_{\nu+1})-\varphi(T_\nu)\right],\ d\in (\varphi(T_\nu),\varphi(T_{\nu+1})),
\end{eqnarray*}
i.e. (see \cite{6}, (5.2); $1.9T<\varphi(T)<2T$)
\be \label{form4.4}
\varphi(T_{\nu+1})-\varphi(T_\nu)=\mcal{O}\left(\frac{1}{\ln T\nu}\right) .
\ee
Next we have
\begin{eqnarray*}
& &
F[\varphi(T_{\nu+1})]-F[\varphi(T_\nu)]=\left( \frac{1}{2}\ln\frac{\varphi(T_\nu)}{2}-\frac{a}{2}\right)
\left[ \varphi(T_{\nu+1})-\varphi(T_\nu)\right]+ \\
& &
\mcal{O}\left\{ \frac{(\varphi(T_{\nu+1})-\varphi(T_\nu))^2}{T}\right\}+\mcal{O}\left(\frac{\ln T_\nu}{T_\nu}\right)=\omega
\end{eqnarray*}
by (\ref{form4.2}), (\ref{form4.3}), i.e. (see (\ref{form4.4}))
\bdis
\frac{1}{2}[\varphi(T_{\nu+1})-\varphi(T_\nu)]\left( \ln\frac{\varphi(T_\nu)}{2}-a\right)=\omega+\mcal{O}\left(\frac{\ln T_\nu}{T_\nu}\right)+
\mcal{O}\left(\frac{1}{T_\nu\ln^2T_\nu}\right),
\edis
from which the formula (\ref{form2.2}) follows.

\begin{remark}

Proof of Theorem 2 is similar to the proof of Theorem 1.

\end{remark}

\thanks{I would like to thank Michal Demetrian for helping me with the electronic version of this work.}

\end{document}